\documentclass[a4paper,12pt, times]{article}
\usepackage[utf8]{inputenc}

\usepackage{mathrsfs}
\usepackage{graphicx}
\usepackage{enumerate}
\usepackage{multicol}
\usepackage{color}
\usepackage{amsmath,amssymb,amscd}
\usepackage{amsthm}

\usepackage{times}

\newtheorem{theorem}{\noindent Theorem}[section]

\newtheorem{lemma}{\noindent Lemma}[section]
\newtheorem{corollary}{\noindent Corollary}[section]

\newtheorem{remark}{\indent Remark}[section]

\def\su{\subset}
\def\stb{,\ldots ,}
\def\cd{\cdot}
\def\se{\setminus}
\def\emp{\emptyset}
\def\del{\partial}
\def\cl{{\rm cl}\, }
\def\deg{{\rm deg}\, }
\def\sumin{\sum_{i=1}^n}
\def\sumjn{\sum_{j=1}^n}
\def\ol{\overline}
\def\akkor{\Longrightarrow}
\def\am{^{-1}}
\def\msk{\medskip}

\def\noi{\noindent}

\def\al{\alpha}
\def\De{\Delta}
\def\la{\lambda}
\def\La{\Lambda}

\def\qq{{\mathbb Q}}
\def\zz{{\mathbb Z}}
\def\cc{{\mathbb C}}

\def\aa{{\cal A}}
\def\dd{{\cal D}}
\def\oo{{\cal O}}

\begin{document}
\setcounter{page}{1}
\newcommand\balline{\small First author  and Second author}
\newcommand\jobbline{\small Short title}

\vspace{-4cm} 

\vspace{.4cm}

\title{Derivations and differential operators on rings and fields}

\author{Gergely Kiss and Mikl\'os Laczkovich}
\date{}

\maketitle

\begin{abstract} {\small
Let $R$ be an integral domain of characteristic zero. We prove that
a function $D\colon R\to R$ is a derivation of order $n$ if and only
if $D$ belongs to the closure of the set of differential operators
of degree $n$ in the product topology of $R^R$, where the image
space is endowed with the discrete topology. In other words, $f$ is
a derivation of order $n$ if and only if, for every finite set $F\su
R$, there is a differential operator $D$ of degree $n$ such that
$f=D$ on $F$. We also prove that if $d_1 \stb d_n$ are nonzero
derivations on $R$, then $d_1 \circ \ldots \circ d_n$ is a
derivation of exact order $n$.}
\end{abstract}



\section{Introduction and main results}
By a ring we mean a commutative ring with unit. An integral domain
is a ring with no zero-divisors other than 0. The ring $R$ has
characteristic zero if $n\cd x\ne 0$ for every $x\in R\se \{ 0\}$
and for every positive integer $n$.

A {\it derivation} on a ring $R$ is a map $d:R\to R$ such that
\begin{equation}\label{eder1}
d(x+y)=d(x)+d(y) \qquad \text{and} \qquad d(xy)=d(x)y+d(y) x
\end{equation}
for every $x,y\in R$. Derivations of higher order are defined by
induction as follows (cf. \cite{UR}).

Let $R$ be a ring. The identically $0$ function defined on $R$ is
called the derivation of order $0$. Let $n>0$, and suppose we have
defined the derivations of order at most $n-1$. A function $D:R\to
R$ is called a {\it derivation of order at most $n$}, if $D$ is
additive and satisfies
\begin{equation}\label{edern}
    D(xy)-D(x)y-D(y)x=B(x,y)
\end{equation}
for every $x,y\in R$, where $B(x,y)$ is a derivation of order at
most $n-1$ in each of its variables. We denote by $\dd ^n (R)$ the
set of derivations of order at most $n$ defined on $R$. We may write
$\dd ^n$ instead of $\dd ^n (R)$ if the ring $R$ is clear from the
context. We say that the order of a derivation $D$ is $n$ if $D\in
\dd ^n \se \dd ^{n-1}$. (We have $\dd ^{-1}=\emp$ by definition).

Clearly, a function $d\colon R\to R$ is a derivation if and only if
$d\in \dd _1$.

Now we define differential operators on a ring $R$. We say that the
map $D:R\to R$ is a {\it differential operator of degree at most
$n$} if $D$ is the linear combination, with coefficients from $R$,
of finitely many maps of the form $d_1 \circ  \ldots \circ  d_k$,
where $d_1 \stb d_k$ are derivations on $R$ and $k\le n$. If $k=0$
then we interpret $d_1 \circ  \ldots \circ  d_k$ as the identity
function on $R$. We denote by $\oo ^n (R)$ the set of differential
operators of degree at most $n$ defined on $R$. We may write $\oo
^n$ instead of $\oo ^n (R)$ if the ring $R$ is clear from the
context. We say that the degree of a differential operator $D$ is
$n$ if $D\in \oo ^n \se \oo ^{n-1}$ (where $\oo ^{-1}=\emp$ by
definition).

The term ``differential operator'' is justified by the following
fact. Let $K=\qq (t_1 \stb t_k  )$, where $t_1 \stb t_k$ are
algebraically independent over $\qq$. Then $K$ is the field of all
rational functions of $t_1 \stb t_k$ with rational coefficients. It
is clear that $d_i =\tfrac{\del}{\del t_i}$ is a derivation on $K$
for every $i=1\stb k$. Therefore, every differential operator
\begin{equation}\label{euj1}
D=\sum_{i_1 +\ldots + i_k \le n} c_{i_1 \stb i_k} \cd
\frac{\partial^{i_1+\dots+i_k}}{\partial t_1^{i_1}\cdots \partial
t_k^{i_k}},
\end{equation}
where the coefficients $c_{i_1 \stb i_k}$ belong to $K$, is a
differential operator of degree at most $n$. The converse is also
true: if $D$ is a differential operator of degree at most $n$ on the
field $K=\qq (t_1 \stb t_k  )$, then $D$ is of the form \eqref{euj1}
(see \cite[Proposition 3.2]{KL} and the proof of Lemma \ref{l3}
below).

\begin{remark}\label{r1}
{\rm If $d$ is a derivation on $R$, then $c\cd d$ is also a
derivation for every $c\in R$. Thus every differential operator is
the sum of terms of the form $d_1 \circ  \ldots \circ  d_k$, where
$k\ge 1$ and $d_1 \stb d_k$ are derivations, and of a term $c\cd j$,
where $c\in R$ and $j$ is the identity function. Since $d(1)=0$ for
every derivation $d$, it follows that a differential operator $D$
satisfies $D(1)=0$ if and only if the term $c\cd j$ is missing; that
is, if $D$ is the sum of terms of the form $d_1 \circ  \ldots \circ
d_k$, where $k\ge 1$ and $d_1 \stb d_k$ are derivations. We denote
by $\oo ^n_0$ the set of all differential operators $D$ of degree at
most $n$ satisfying $D(1)=0$. }
\end{remark}

Let $G$ be an Abelian semigroup, and let $H$ be an Abelian group.
The {\it difference operator} $\De _g$ $(g\in G)$ is defined by $\De
_g f(x)=f(x+g)- f(x)$ for every $f\colon G \to H$ and $x\in G$. A
function $f:G\to H$ is a {\it generalized polynomial}, if there is a
$k$ such that $\De _{g_1}\ldots \De _{g_{k+1}}f=0$ for every $g_1
\stb g_{k+1} \in G$. The smallest $k$ for which this holds for every
$g_1 \stb g_{k+1} \in G$ is the {\it degree} of the generalized
polynomial $f$, denoted by $\deg f$. The degree of the identically
zero function is $-1$ by definition. It is clear that the nonzero
constant functions are generalized polynomials of degree $0$, and
the nonconstant additive functions; that is, the nonzero
homomorphism from $G$ to $H$, are generalized polynomials of degree
$1$.

If $X,Y$ are nonempty sets, then $Y^X$ denotes the set of all maps
$f\colon X\to Y$. We endow the space $Y$ with the discrete topology,
and $Y^X$ with the product topology. The closure of a set $\aa \su
Y^X$ with respect to the product topology is denoted by $\cl \aa$.
Clearly, a function $f\colon X\to Y$ belongs to $\cl \aa$ if and
only if, for every finite set $F\su X$ there is a function $g\in
\aa$ such that $f(x)=g(x)$ for every $x\in F$.

It is clear that a function $f\colon G \to H$ is a generalized
polynomial of degree at most $n$ if and only if, for every finite
set $F\su G$, there is a generalized polynomial $h$ of degree at
most $n$ such that $f=h$ on $F$. This means that {\it the set of
generalized polynomials of degree at most $n$ is closed in $H^G$.}

If $R$ is a ring, then we denote by $R^*$ the Abelian semigroup
$R\se \{ 0\}$ under multiplication. We denote by $j$ the identity
function on $R$.

In this note our aim is to prove that, for every integral domain of
characteristic zero and for every positive integer $n$, we have $\dd
^n =\cl \oo ^n_0$. That is, {\it a map $D\colon R\to R$ is a
derivation of order at most $n$ if and only if $D$ belongs to the
closure of the set of all differential operators of degree at most
$n$ satisfying $D(1)=0$.} More precisely, we prove the following
result.

\begin{theorem}\label{t1}
Let $R$ be an integral domain of characteristic zero, $K$ its field
of fractions, and let $n$ be a positive integer. Then, for every
function $D\colon R\to R$, the following are equivalent.
\begin{enumerate}[{\rm (i)}]
\item $D\in \dd ^n (R)$.
\item $D\in \cl (\oo ^n_0 (R))$.
\item $D$ is additive on $R$, $D(1)=0$, and $D/j$, as a map from the semigroup
$R^*$ to $K$, is a generalized polynomial of degree at most $n$.
\end{enumerate}
\end{theorem}

As an immediate consequence of the theorem above we find the
following corollary.

\begin{corollary}\label{c1}
Let $R$ be an integral domain of characteristic zero, $K$ its field
of fractions, and let $n$ be a positive integer. Then, for every
function $D\colon R\to R$, the following are equivalent.
\begin{enumerate}[{\rm (i)}]
\item $D\in \dd ^n (R)\se \dd^{n-1} (R)$.
\item $D\in (\cl \oo ^n_0 (R)) \se \cl (\oo ^{n-1}_0 (R))$.
\item $D$ is additive on $R$, $D(1)=0$, and $D/j$, as a map from the semigroup
$R^*$ to $K$, is a generalized polynomial of degree $n$.
\end{enumerate}
\end{corollary}

Indeed, suppose $D\in \dd ^n \se \dd^{n-1}$. Then, by Theorem
\ref{t1}, we have $D\in \cl \oo ^n_0$. If $D\notin \cl (\oo ^n_0 )
\se \cl (\oo ^{n-1}_0 )$, then $D\in \cl \oo ^{n-1}_0$. This implies
$D\in \dd ^{n-1}$, which is impossible. Therefore, (i) of Corollary
\ref{c1} implies (ii) of Corollary \ref{c1}. The other implications
can be shown similarly.

\begin{remark}\label{r2}
{\rm Theorem \ref{t1} and Corollary \ref{c1} do not hold without
assuming
  that $R$ is of characteristic zero. Consider the following example.

Let $F_2$ denote the field having two elements, and let $R=F_2 [x]$
be the ring of polynomials with coefficients from $F_2$. We put
$$D\left( \sum_{i=0}^n a_i \cd x^i \right) =\sum_{i=2}^n \frac{i(i-1)}{2} \cd
a_i \cd x^{i-2}$$ for every $n\ge 0$ and $a_0 \stb a_n \in F_2$. It
is easy to check that $D$ is a derivation of order at most two on
$R$. Since $D(x)=0$ and $D(x^2)=1$, it follows that $D$ is not a
derivation, and thus $D\in \dd ^2 \se \dd^1$.

On the other hand, if $d_1$ and $d_2$ are arbitrary derivations on
$R$, then $d_1 \circ d_2$ is also a derivation. Indeed,
$$d_1 (d_2 (x^k ))=d_1 (k\cd x^{k-1} \cd d_2 (x))=k(k-1)\cd x^{k-2} \cd d_1 (x)\cd d_2 (x)+ k\cd x^{k-1} \cd d_1 (d_2 (x)) $$
for every $k\ge 2$. Since $k(k-1)$ is even, we find that
\begin{equation}\label{e5}
(d_1 \circ d_2 ) (x^k )= k\cd x^{k-1}\cd a
\end{equation}
for every $k\ge 2$, where $a=d_1 (d_2 (x))\in R$. It is easy to
check that \eqref{e5} is true for $k=0$ and $k=1$ as well. Since
derivations are additive, \eqref{e5} gives $d_1 (d_2 (p))=a\cd
\tfrac{\del p}{\del x}$ for every $p\in R$, and thus $d_1 \circ d_2
\in \oo ^1_0$. This implies that $\oo ^2_0 =\oo ^1_0$, and thus $\dd
^2$ is strictly larger than $\oo ^2_0$. }
\end{remark}

\begin{remark}\label{r4}
{\rm In the proof of Theorem \ref{t1} the crucial step is to show
that if $R$ is of characteristic zero and the transcendence degree
of the field of fractions $K$ of $R$ over $\qq$ is finite, then $\dd
^n =\oo ^n _0$ (see Lemma \ref{l5}). Comparing to Theorem \ref{t1}
we find that under these conditions, for every function $f\colon
R\to R$ we have

$(f\in \dd^n \se \dd ^{n-1})\iff (f\in \oo^n_0 \se \oo^{n-1}_0 )
\iff D$ is additive on $K$, $D(1)=0$, and $D/j$, defined on the
group $K^*$, is a generalized polynomial of degree $n$.}
\end{remark}

We also prove that for every integral domain $R$ of characteristic
zero, if there are nonzero derivation on $R$, then the sets $\dd ^n
\se \dd^{n-1}$ are nonempty; that is, there are derivations of any
given order. More precisely, we prove the following \footnote{Added
in proof: it came to our notice recently that the statement of
Theorem 1.2 was also proved, using different methods, by Bruce
Ebanks in his submitted paper "Derivations and Leibniz differences
on rings".}.

\begin{theorem}\label{t2}
Let $R$ be an integral domain of characteristic zero, and let $n$ be
a positive integer. If $d_1 \stb d_n$ are nonzero derivations on
$K$, then $d_1 \circ \ldots \circ d_n \in \dd ^n \se \dd^{n-1}$.
\end{theorem}

(For integral domains of characteristic zero this generalizes
\cite[Remark 3]{GKV}, where the case $d_1 =\ldots =d_n$ is
considered.)

\begin{remark}\label{r3}
{\rm The statement of the theorem above does not hold without
assuming that $R$ is of characteristic zero. Consider the example
described in Remark \ref{r2}. Clearly, $d(p)=\tfrac{\del p}{\del x}$
($p\in R$) defines a nonzero derivation on $R$. However, as we saw
in Remark \ref{r2}, $d\circ d$ is a derivation of order $1$.

The statement of the theorem is not true for rings in general; not
even for rings of characteristic zero. Let $R=\qq [x]\times \qq
[x]$, and put $d_1 (p,q)= (\tfrac{\del p}{\del x} ,0)$ and $d_2
(p,q)= (0, \tfrac{\del q}{\del x})$ for every $(p,q)\in R$. Then
$d_1$ and $d_2$ are nonzero derivations on $R$, but $d_1 \circ d_2
=0$. }
\end{remark}

\section{Lemmas}
\begin{lemma}\label{l1}
For every ring $R$ and for every nonnegative integer $n$, the set
$\dd ^n$ is closed in $R^R$.
\end{lemma}

\proof We prove by induction on $n$. If $n=0$, then $\dd ^0 =\{ 0\}$
is closed. Let $n>0$, and suppose that $\dd ^{n-1}$ is closed. Let
$f\in \cl \dd ^n$ be arbitrary. We have to prove that $f\in \dd ^n$;
that is, for every fixed $y\in R$, the map $x\mapsto
g(x)=f(xy)-yf(x)-xf(y)$ belongs to $\dd ^{n-1}$. By the induction
hypothesis, it is enough to show that $g\in \cl \dd ^{n-1}$; that
is, for every finite set $F\su R$ there is a function $h\in \dd
^{n-1}$ such that $g(x)=h(x)$ for every $x\in F$.

If $F$ is finite, then so is $A=F\cup \{ xy\colon x\in F\} \cup \{
y\}$. Since $f\in \cl \dd ^n$, there is a function $D\in \dd ^n$
such that $f(z)=D(z)$ for every $z\in A$. If $x\in F$, then
$x,y,xy\in A$, and thus
$$g(x)=f(xy)-yf(x)-xf(y)=D(xy)-yD(x)-xD(y).$$
The function $x\mapsto h(x)=D(xy)-yD(x)-xD(y)$ belongs to $\dd
^{n-1}$, as $D\in \dd ^n$. Since $g(x)=h(x)$ for every $x\in F$, the
lemma is proved. \hfill $\square$

\begin{lemma}\label{l4}
For every ring $R$ we have $\cl \oo ^n_0 \su \dd ^n$.
\end{lemma}

\proof Since $\dd ^n$ is closed by Lemma \ref{l1}, it is enough to
show that $\oo ^n_0 \su \dd ^n$. Let $D$ be a differential operator
of degree at most $n$ satisfying $D(1)=0$. According to Remark
\ref{r1}, $D$ is the sum of terms of the form $d_1 \circ  \ldots
\circ  d_k$, where $1\le k\le n$ and $d_1 \stb d_k$ are derivations.
Since $\dd ^n$ is a linear space, it is enough to show that $d_1
\circ  \ldots \circ  d_k \in \dd ^k$ whenever $k\ge 1$ and $d_1 \stb
d_k$ are derivations. This, in turn, is easy to prove by induction
on $k$. \hfill $\square$

The statement of the following lemma is probably known. In order to
make these notes as self-contained as possible, we provide the
proof.

\begin{lemma}\label{l7}
Let $G$ be an Abelian semigroup, and let $K$ be a field. If $p\colon
G\to K$ is a generalized polynomial of degree $n\ge 0$ and $a\colon
G\to K$ is a nonzero additive function, then $p\cd a$ is a
generalized polynomial of degree at most $n+1$.

If $K$ is of characteristic zero, then $\deg (p\cd a)=n+1$.
\end{lemma}

\proof We prove by induction on $n$. If $n=0$, then $p$ is a nonzero
constant, and $p\cd a$ is a nonzero additive function, hence a
generalized polynomial of degree $1$.

Let $n>0$, and suppose that the statement is true for $n-1$. Let $p$
be a generalized polynomial of degree $n$. We have
\begin{equation}\label{e3}
\De _g (p \cd a)(x)= a(x)\cd \De _g p (x)+a(g)\cd p (x+g)
\end{equation}
for every $x,g\in G$. Since $\deg  \De _g p (x) \le n-1$, it follows
from the induction hypothesis that $\deg (a(x)\cd \De _g p (x)) \le
n$. Therefore, by \eqref{e3}, we have $\deg \De _g (p\cd a) \le n$
for every $g\in G$, and thus $\deg (p\cd a)\le n+1$. We have to
prove that if $K$ is characteristic zero, then $\deg (p\cd a)\ge
n+1$.

Since the image space $K$ is a torsion free and divisible Abelian
group, it follows from Djokovi\'c's theorem \cite{Dj} that $p=P_n
+\ldots +P_1 +P_0$, where $P_i$ is a monomial of degree $i$ for
every $i=1\stb n$, and $P_0$ is constant. Then there is a symmetric
function $A(x_1 \stb x_n )$, additive in each of its variables, such
that $P_n (x)=A(x\stb x)$ $(x\in G)$. Since $q=p-P_n$ is a
generalized polynomial of degree $\le n-1$, it follows from the
induction hypothesis that $\deg (q\cd a)\le n$. Therefore, in order
to prove $\deg (p\cd a)\ge n+1$, it is enough to show that $\deg
(P_n \cd a)=n+1$.

First we show that there exists an element $g\in G$ such that $P_n
(g)\ne 0$ and $a(g)\ne 0$. By assumption, there is an $x\in G$ such
that $a(x)\ne 0$. Since $\deg P_n =n\ge 0$, it follows that $P_n$ is
nonzero. Let $y\in G$ be such that $P_n (y)\ne 0$. Now $a(kx+y)=k\cd
a(x)+ a(y)$ for every positive integer $k$. Since $a(x),a(y)\in K$
and $a(x)\ne 0$, we have $a(kx+y)\ne 0$ for every $k$ with at most
one exception.

Using the fact that $A(x_1 \stb x_n )$ is symmetric and additive in
each of its variables, we find
\begin{equation}\label{e4}
 P_n (kx+y) = \sum_{i=0}^n \binom{n}{i} A_i (kx,y)
\end{equation}
for every positive integer $k$, where
$$A_i (kx,y) =A (\underbrace{kx \stb kx}_{i} ,\underbrace{y\stb y}_{k-i} ) =
k^i \cd A (\underbrace{x \stb x}_{i} ,\underbrace{y\stb y}_{k-i} )
.$$ Therefore, by \eqref{e4}, $Q(kx+y)$ is a polynomial of $k$ with
coefficients from $K$. Since the constant term of this polynomial is
$A(y\stb y)\ne 0$, $Q(kx+y)$ is not the identically zero polynomial,
and thus $ P_n (kx+y)\ne 0$ for all but finitely many $k$.
Therefore, we may choose a $k$ such that $P_n (g)\ne 0$ and $a(g)\ne
0$, where $g=kx+y$.

Let $Q=P_n \cd a$, and suppose that $\deg Q\le n$. Then $Q=Q_n
+\ldots +Q_1 + Q_0$, where $Q_i$ is a monomial of degree $i$ for
every $i=1\stb n$, and $Q_0$ is constant. For every $i=1\stb n$,
there is there is a symmetric function $B_i (x_1 \stb x_i )$,
additive in each of its variables, such that $Q_i (x)=B_i (x\stb x)$
$(x\in G)$. Then
$$Q(k \cd g)=Q_0 +\sumin B_i (kg\stb kg)=Q_0 +\sumin k^i \cd B_i (g\stb kg)$$
for every positive integer $k$. Therefore, the map $k\mapsto Q(k \cd
g)$ is a polynomial of degree $\le n$ with coefficients from $K$.
However,
$$Q(k \cd g)=k^n \cd A(g\stb g)\cd k\cd a(g)=k^{n+1} \cd A(g\stb g)\cd a(g)$$
is a polynomial of degree $n+1$. This is a contradiction, proving
$\deg Q=n+1$. \hfill $\square$

\begin{lemma}\label{l2}
Let $R$ be an integral domain, and let $K$ be its field of
fractions. If $d_1 \stb d_n$ are nonzero derivations on $R$ and
$D=d_1 \circ \ldots \circ d_n$, then $D/j$, as a map from the
semigroup $R^*$ to $K$, is a generalized polynomial of degree at
most $n$.

If $R$ is of characteristic zero, then $\deg D/j =n$.
\end{lemma}

\proof We prove by induction on $n$. If $n=1$, then $D$ is a nonzero
derivation. It is clear that in this case $D/j$ is additive, hence a
generalized polynomial of degree at most $1$ on the semigroup $R^*$.
Suppose $\deg D/j \le 0$. Then $D/j$ is constant on $R^*$, and thus
$D=c\cd j$ on $R$, where $c\in R$ is a constant. Since $D$ is a
derivation, we have $c=D(1)=0$ and $d=0$, a contradiction. Thus
$\deg D/j=1$.

Suppose that $n>1$, and the statement is true for $n-1$. Let $d_1
\stb d_n$ be nonzero derivations on $R$. By the induction
hypothesis, $(d_2 \circ \ldots \circ d_n )/j =p$ is a generalized
polynomial of degree at most $n-1$. Since $d_1$ is a derivation, we
have
$$D(x)=(d_1 \circ  \ldots \circ  d_n )(x)=d_1 (p(x)\cd x)=d_1 (p(x))\cd x+
p(x)\cd d_1 (x)$$ for every $x\in R^*$. Thus
\begin{equation}\label{e1}
D/j =(d_1 \circ p)+p\cd (d_1 /j)
\end{equation}
on $R^*$. Since $p\colon R^* \to K$ is a generalized polynomial of
degree $\le n-1$ and $d_1 :R\to R$ is additive, it follows that $d_1
\circ p$ is a generalized polynomial of degree $\le n-1$ on $R^*$.
(This is because, if $G$ is an Abelian semigroup, $H$ is an Abelian
group, $p\colon G\to H$ is a generalized polynomial of degree $k$,
and $d\colon H\to H$ is additive, then $d\circ p$ is a generalized
polynomial of degree at most $k$.)

If $R$ is of characteristic zero, then so is $K$. In this case $p\cd
(d_1 /j)$ is a generalized polynomial of degree $n$ by Lemma
\ref{l7}, since $d_1 /j$ is nonzero and additive on $R^*$.
Therefore, $D/j$ is a generalized polynomial of degree $n$. \hfill
$\square$

\begin{lemma}\label{l6}
Let $R$ be an integral domain, and let $K$ be its field of
fractions. If $D\in \cl \oo ^n_0 (R)$, then $D/j$, as a map from the
semigroup $R^*$ to $K$, is a generalized polynomial of degree at
most $n$.
\end{lemma}

\proof Let $D\in \cl \oo ^n_0$ be given. As the set of generalized
polynomials of degree $\le n$ is closed, it is enough to show that
for every finite set $F\su R^*$ there is a generalized polynomial
$p\colon R^* \to K$ such that $\deg p \le n$ and $D/j=p$ on $F$.
Since $D\in \cl \oo ^n_0$, there is an $f\in \oo ^n_0$ such that
$D=f$ on $F$. It is clear from Remark \ref{r1} and Lemma \ref{l2}
that $f/j$ is a generalized polynomial of degree at most $n$. Now we
have $D/j=f/j$ on $F$, completing the proof. \hfill $\square$

The statement of the following lemma is proved, in a different
context, in Lemma 3.3 of \cite{KL}. We give the proof adjusted to
our purposes.

\begin{lemma}\label{l3}
Let $R$ be a subring of $\cc$, let $K\su \cc$ be its field of
fractions, and suppose that the transcendence degree of $K$ over
$\qq$ is finite. Let the map $D:R\to R$ be additive. If $D/j$, as a
map from the semigroup $R^*$ to $\cc$ is a generalized polynomial of
degree at most $n$, then $D\in \oo ^n$.
\end{lemma}

\proof Let $k$ be the transcendence degree of $K$ over $\qq$, and
let the elements $u_1 \stb u_k \in K$ be algebraically independent
over $\qq$. Let $u_i =a_i /b_i$, where $a_i ,b_i \in R$ for every
$i=1\stb k$. Then the field $\qq (a_1 ,b_1 \stb a_k ,b_k )$ has
transcendence degree $k$ over $\qq$, and thus we can chose elements
$t_1 \stb t_k \in \{ a_1 ,b_1 \stb a_k ,b_k \} \su R^*$ such that
$t_1 \stb t_k$ are algebraically independent over $\qq$.

By assumption, the function $p=D/j$ is a generalized polynomial of
degree $\le n$ on $R^*$. By Djokovi\'c's theorem, we have $p=P_n
+\ldots +P_1 + P_0$, where $P_j$ is a monomial of degree $j$ for
every $j=1\stb n$, and $P_0$ is constant. Using the fact that $P_j
(x)=A_j (x\stb x)$, where $A_j (x_1 \stb x_j )$ is symmetric and
additive in each of its variables, it is easy to see that for every
$j=1\stb n$ there is a homogeneous polynomial $\ol{p}_j \in K [x_1
\stb x_k ]$ of degree $j$ such that
$$P_j  \left( t_1^{i_1} \cdots t_k^{i_k} \right) =\ol{p} _j (i_1 \stb i_k)$$
whenever $i_1 \stb i_k$ are nonnegative integers. (Note that the
semigroup operation in $R^*$ is multiplication.) Putting $\ol p =P_0
+\sumjn \ol{p}_j$ we find that $\ol p \in K[x_1 \stb x_k ]$, and
\begin{equation*}\label{e9}
\ol p \left( t_1^{i_1} \cdots t_k^{i_k} \right) =q (i_1 \stb i_k)
\end{equation*}
for every $i_1 \stb i_k \ge 0$. We shall use the notation
$x^{[0]}=1$ and $x^{[j]}=x(x-1)\cdots (x-j+1)$ for every
$j=1,2,\ldots $ and $x\in \zz .$ It is easy to see that every
polynomial belonging to $K [x_1 \stb x_k ]$ and of degree $\le n$
can be written in the form $\sum c_j \cd x_1^{[j_1 ]} \cdots
x_k^{[j_k ]}$, where $j=(j_1 \stb j_k )$ runs through the set of
$k$-tuples of nonnegative integers with $j_1 +\ldots +j_k \le n$,
and in each term the coefficient $c_j$ belongs to $K$. Therefore,
the polynomial $\ol p$ also has such a representation. Then we have
\begin{equation}\label{e8}
\begin{split}
D\left( t_1^{i_1} \cdots t_k^{i_k} \right) & =
p\left( t_1^{i_1} \cdots t_k^{i_k} \right) \cd  t_1^{i_1} \cdots t_k^{i_k}  = \\
& =\sum c_j \cd i_1^{[j_1 ]} \cdots i_k^{[j_k ]} \cd t_1^{i_1} \cdots t_k^{i_k} = \\
& = \sum c_j \cd  t_1^{j_1} \cdots t_k^{j_k} \cd i_1^{[j_1 ]} \cdots
i_k^{[j_k ]} \cd t_1^{i_1 -j_1 } \cdots t_k^{i_k -j_k} =\\
&=E \left( t_1^{i_1} \cdots t_k^{i_k} \right)
\end{split}
\end{equation}
for every $i_1 \stb i_k \ge 0$, where $E$ is the differential
operator
$$\sum c_j \cd t_1^{j_1} \cdots t_k^{j_k}  \cd
\frac{\partial^{j_1+\dots+j_k}}{\partial t_1^{j_1}\cdots \partial
t_k^{j_k}}.$$ By extending the derivations $\del /\del t_i$ to $K$,
we can extend $E$ to $K$ as a differential operator $\overline{E}$
of degree at most $n$. Then $\ol E$ is additive on $K$, and
$\overline{E}/j$ is a generalized polynomial on $K^*$ by Lemma
\ref{l2}. Let $q(0)=0$, and let $q(x)=p(x) -\ol E (x)/x$ for every
$x\in R^*$. Then $q\cd j =D- \ol E$ is additive on $R$, and $q$ is a
generalized polynomial on $R^*$. Let $G$ denote the semigroup
generated by the elements $t_1 \stb t_k$. Then $q$ vanishes on $G$
by \eqref{e8}. From these conditions it follows that $q= 0$ on $R$.
This is proved in \cite[Lemma 3.6]{KL} under the stronger condition
that $G$ is the group (and not the semigroup) generated by $t_1 \stb
t_k$. One can see that the same argument works in our more general
case as well; however, for the sake of completeness we give the
proof in the appendix. Thus we have $q=0$; that is, $D=\ol E$ on
$R$, which completes the proof. \hfill $\square$

\begin{lemma}\label{l5}
Let $R$ be a subring of $\cc$, let $K\su \cc$ be its field of
fractions, and suppose that the transcendence degree of $K$ over
$\qq$ is finite. Then $\dd ^n (R) =\oo ^n_0 (R)$.
\end{lemma}

\proof By Lemma \ref{l4}, we only have to show that $\dd ^n \su \oo
^n_0$. It is easy to prove, by induction on $n$ that if $D\in \dd
^n$, then $D(1)=0$. Therefore, it is enough to show that if $D\in
\dd ^n$, then $D$ is a differential operator of degree at most $n$.
We prove by induction on $n$.

The statement is obvious if $n=0$. Let $n>0$, and suppose that the
statement is true for $n-1$. Let $D$ be a derivation of order at
most $n$. By Lemma \ref{l3}, it is enough to show that $p=D/j$,
defined on the semigroup $R^*$, is a generalized polynomial of
degree at most $n$. Let $y\in R^*$ be fixed. Dividing \eqref{edern}
by $xy$ we obtain
$$\frac{D(xy)}{xy}-\frac{D(x)}{x}-\frac{D(y)}{y}=\frac{B(x,y)}{xy},$$
and thus $p(xy)-p(x)-p(y)=B(x,y)/{xy}$ for every $x\in K^*$.
Therefore we have
\begin{equation}\label{e2}
\De _y p (x)=p(y)+\frac{1}{y}\cd \frac{B(x,y)}{x}
\end{equation}
on $R^*$. The map $x\mapsto B(x,y)$ is a derivation of order at most
$n-1$. We also have $B(1,y)=0$ by $D(1)=0$. Therefore, by Lemma
\ref{l2}, the map $x\mapsto B(x,y)/x$ is a generalized polynomial of
degree at most $n$. Then so is $\De _y p$ by \eqref{e2}. Since this
is true for every $y\in K^*$, it follows that $p$ is a generalized
polynomial of degree at most $n$. \hfill $\square$

\section{Proof of Theorems \ref{t1} and \ref{t2}.}
First we prove Theorem \ref{t1}. The implication (ii)$\akkor$(iii)
is proved in Lemma \ref{l6}.

\msk \noi (iii)$\akkor$(ii): Suppose that $D$ is additive, $D(1)=0$,
and $D/j$ is a generalized polynomial of degree at most $n$. In
order to prove $D\in \cl \oo ^n_0$, we have to show that for every
finite set $F\su K$ there is a function $f\in \oo ^n_0$ such that
$D=f$ on $F$. Let $F\su K$ be finite, and let $L$ denote the
subfield of $K$ generated by $F$. Obviously, the transcendence
degree of $L$ over $\qq$ is finite. It is well-known that every
field of characteristic zero and having finite transcendence degree
over $\qq$ is isomorphic to a subfield of $\cc$. Therefore, we may
assume that $L\su \cc$. Thus, by Lemma \ref{l3}, the restriction
$D|_L$ of $D$ to the field $L$ is a derivation of order at most $n$.
Since $D(1)=0$, we also have $D|_L \in \oo ^n _0 (L)$. It is
well-known that every derivation on $L$ can be extended to $K$ as a
derivation (see \cite[pp. 351-352]{K}). This implies that every
differential operator on $L$ of degree at most $n$ can be extended
to $K$ as a differential operator of degree at most $n$. If $f$ is
such an extension of $D|_L$, then, obviously, $D(x)=f(x)$ for every
$x\in F$. This proves (iii)$\akkor$(ii).

\msk \noi (ii)$\akkor$(i): This is Lemma \ref{l4}.

\msk \noi (i)$\akkor$(ii): Let $D\in \dd ^n$. In order to prove
$f\in \cl \oo ^n_0$ we have to show that for every finite set  $F\su
K$ there is a function $f\in \oo ^n_0$ such that $D=f$ on $F$. Let
$L$ denote the field generated by $F$. Obviously, the transcendence
degree of $L$ over $\qq$ is finite. Thus, by Lemma \ref{l5}, the
restriction $D|_L$ of $D$ to the field $L$ is a derivation of order
at most $n$, vanishing at $1$. Let $f$ be an extension of $D|_L$ to
$K$ as a function $f\in \oo ^n_0$. Then, obviously, $D(x)=f(x)$ for
every $x\in F$. This proves (i)$\akkor$(ii). \hfill $\square$

The statement of Theorem \ref{t2} is an immediate consequence of
Corollary \ref{c1} and Lemma \ref{l2}. \hfill $\square$

\section{Appendix}
\begin{lemma}\label{l8}
Let $R$ be a subring of $\cc$, and let $K\su \cc$ be its field of
fractions. Suppose that the transcendence degree of $K$ over $\qq$
is $k<\infty$, and let the elements $t_1 \stb t_k \in R$ be
algebraically independent over $\qq$. Let $f\colon R\to \cc$ be
additive on $R$ (with respect to addition) and such that $q=f/j$, as
a map from the semigroup $R^*$ to $\cc$ is a generalized polynomial.
If $f=0$ on the semigroup $G$ generated by $t_1 \stb t_k$, then
$f=0$ on $R$.
\end{lemma}

\proof We prove by induction on $\deg q$. If $\deg q=0$, then $q$ is
constant. Since $f= 0$ on $G$, we have $q=0$ on $G$, and thus $q=0$
on $R$.

Suppose $m=\deg q>0$, and that the statement is true for degrees
less than $m$. Let $g\in G$ be fixed, and put $f_1 (x)=g\am
f(gx)-f(x)$ $(x\in R)$. Then $f_1$ is additive on $R$. Also, $f_1
/j$ is a generalized polynomial on $R^*$, since
\begin{equation*}\label{e6}
  \frac{f_1 (x)}{x} =\frac{g\am f(gx)-f(x) }{x } =\frac{f(gx)}{gx}-
  \frac{f(x)}{x}=q(gx)-q(x) =\De _g q(x)
\end{equation*}
for every $x\in R^*$. Since  $\deg (f_1 /j)=\deg \De_g q \le m-1$
and $f_1 = 0$ on $G$, it follows from the induction hypothesis that
$f_1 = 0$ on $R$. Thus $f(gx)=g\cd f(x)$ for every $g\in G$ and
$x\in R$. By the additivity of $f$ we obtain
\begin{equation}\label{e7}
f(cx)=c\cd f(x) \qquad  (c\in \qq [t_1 \stb t_k ] , \ x\in R).
\end{equation}
Since the transcendence degree of $K$ over $\qq$ is $k$ and $t_1
\stb t_k$ are algebraically independent over $\qq$, it follows that
every element of $K$ is algebraic over $\qq (t_1 \stb t_k )$. Let
$\al \in R$ be arbitrary. Then $\al$ is algebraic over the field
$\qq (t_1 \stb t_k )$, and there are elements $c_0 \stb c_N \in \qq
[t_1 \stb t_k ]$ such that
\begin{equation}\label{e10}
c_N \al ^N +\ldots +c_1 \al +c_0 =0,
\end{equation}
where $c_N \ne 0$ and $N$ is minimal. Let $f(\al ^i )=a_i \
(i=0,1,\ldots )$. Multiplying \eqref{e10} by $\al ^{n-N}$ for every
$n\ge N$ we obtain
$$c_N \al ^n +\ldots +c_1\al ^{n-N+1}+c_0 \al ^{n-N}=0.$$
By \eqref{e7} and by the additivity of $f$, this implies
$$c_N a_n +\ldots +c_1 a_{n-N+1}+c_0 a_{n-N}=0$$
for every $n\ge N$. Therefore, the sequence $(a_n )$ satisfies a
linear recurrence relation. It is well-known that $a_n$ can be
uniquely represented in the form $a_n =\sum_{\la \in \La} p_\la
(n)\cd \la ^n$, where $\la$ runs through $\La$, the set of roots of
the characteristic polynomial $\chi (x)=c_N x^N +\ldots +c_0$, and
for every root $\la \in \La$, $p_\la \in \cc [x]$ is a polynomial of
the degree less than the multiplicity of $\la$.

Since $N$ is minimal, the polynomial $\chi$ is irreducible over $\qq
(t_1 \stb t_k )$. Therefore, every $\la$ is a simple root of $\chi$,
and thus
\begin{equation}\label{e11}
a_n =\sum_{\la \in \La} d_\la \cd \la ^n
\end{equation}
for every $n$, where $d_\la$ is a constant for every $\la \in \La$.

Since $q$ is a generalized polynomial on $R^*$ it follows that the
map $n\mapsto q(\al ^n )$ is a polynomial on $\{ 0,1,\ldots \}$.
Now, we have $a_n =f(\al ^n )=q(\al ^n )\cd \al ^n$ for every $n$.
The uniqueness of the representation \eqref{e11} implies that $\al
\in \La$, and the function $n\mapsto q(\al ^n )\ (n=0,1, \ldots )$
is constant. Since $q(1)=f(1)=0$ by $1\in G$, it follows that $q(\al
^n )=0$ for every $n$. In particular, $q(\al )=0$ and $f(\al )=0.$
Since this is true for every $\al \in R$, we obtain $f= 0$ on $R$.
\hfill $\square$

\section*{Acknowledgement}
The authors were supported by the Hungarian National Foundation for
Scientific Research, Grant No. K124749 The first author was
supported by the internal research project R-AGR-0500 of the
University of Luxembourg.

\vspace{4cm}

\noindent\textbf{Gergely Kiss}\\
University of Luxemburg\\
Belval\\
Luxemburg\\
{\tt kigergo57@gmail.com}\\

\noindent\textbf{Mikl\'os Laczkovich}\\
E\"otv\"os Lor\'and University\\
BUdapest\\
Hungary\\
{\tt laczk@cs.elte.hu}

\end{document}